\documentclass[10pt,twoside,reqno]{amsart}
\usepackage{amsmath, amsthm, amscd, amsfonts, amssymb, graphicx, color}
\usepackage[bookmarksnumbered, plainpages, backref]{hyperref}
\numberwithin{equation}{section}
\makeatletter

\renewcommand{\@secnumfont}{\bfseries}

\renewcommand{\section}{\@startsection{section}{1}%
  {0mm}{.7\linespacing\@plus\linespacing}{.5\linespacing}
  {\normalfont\bfseries\centering}}

\newcommand{\bibsection}{\@startsection{section}{1}%
  {0mm}{.7\linespacing\@plus\linespacing}{.5\linespacing}
  {\normalfont\scshape\centering}}

\renewcommand{\@biblabel}[1]{#1.}
\newtheorem{theorem}{Theorem}[section]

\begin{document}

\title{Symmetric identities of higher-order degenerate $q$-Bernoulli polynomials }

\author{Taekyun Kim}
\address{Department of Mathematics, College of Science, Tianjin Polytechnic University, Tianjin City, 300387, China \\ Department of Mathematics, Kwangwoon University, Seoul 139-701, S. Korea}
\email{tkkim@kw.ac.kr}

\author{Hyuck-In Kwon}
\address{Department of Mathematics, Kwangwoon University, Seoul 139-701, Republic
of Korea}
\email{sura@kw.ac.kr}

\begin{abstract}
The purpose of this paper is to give symmetric identities for higher-order degenerate $q$-Bernoulli polynomials arising from the $p$-adic $q$-integral on $\mathbb{Z}_p$.

\end{abstract}

\maketitle
\global\long\def\acl#1#2{\left\langle \left.#1\right|#2\right\rangle }

\global\long\def\acr#1#2{\left\langle #1\left|#2\right.\right\rangle }

\global\long\def\Li{\mathrm{Li}}

\global\long\def\Zp{\mathbb{Z}_{p}}

\markboth{\centerline{\scriptsize Symmetric identities of higher-order degenerate $q$-Bernoulli polynomials}}{\centerline{\scriptsize Symmetric identities of higher-order degenerate $q$-Bernoulli polynomials}}

\section{Introduction}

Let $p$ be a fixed odd prime number. Throughout this paper, $\mathbb{Z}_{p}$, $\mathbb{Q}_{p}$ and $\mathbb{C}_{p}$ will denote the ring of $p$-adic integers, the field of $p$-adic rational numbers and the completion of the algebraic closure of $\mathbb{Q}_{p}$.
The $p$-adic norm $|\cdot|_p$ is normalized as $|p|_p = \frac{1}{p}$. The $q$-analogue of the number $x$ is defined as $\left[x\right]_{q}=\frac{1-q^{x}}{1-q}$. Note that $\lim_{q \rightarrow 1} [x]_q = x.$ Let $q$ be an indeterminate such that $|1-q|_p < p^{- \frac{1}{p-1}}$. As is well known, the Bernoulli polynomials are defined by the generating function to be
\begin{equation}\begin{split}\label{01}
\frac{t}{e^t-1} e^{xt} = \sum_{n=0}^\infty B_n(x) \frac{t^n}{n!}\quad (\textnormal{see} \,\, [7,9,10]).
\end{split}\end{equation}
When $x=0$, $B_n = B_n(0)$ are called ordinary Bernoulli numbers. From \eqref{01}, we can derive the following recurrence relation related to $B_n$ as follows:
\begin{equation}\begin{split}\label{02}
B_0 =1, \quad (B+1)^n - B_n = \begin{cases} 1,&\text{if}\,\, n=1\\0,&\text{if}\,\,n>1, \end{cases}
\end{split}\end{equation}	
with the usual convention about replacing $B^n$ by $B_n$. Let $UD(\mathbb{Z}_p)$ be the space of all uniformly differentiable function on $\mathbb{Z}_p$. For $f \in UD(\mathbb{Z}_p)$, the $p$-adic invariant integral on $\mathbb{Z}_p$ is given by
\begin{equation}\begin{split}\label{03}
I_0(f) = \int_{\mathbb{Z}_p}   f(x) d\mu_{0} (x) = \lim_{N \rightarrow \infty} \frac{1}{p^N} \sum_{x=0}^{p^N-1} f(x),\quad (\textnormal{see} \,\, [8,10]).
\end{split}\end{equation}
From \eqref{03}, we note that
\begin{equation}\begin{split}\label{04}
I_0(f_1) - I_0(f) = f'(0),\quad (\textnormal{see} \,\, [8,10]),
\end{split}\end{equation}
where $f_1(x) = f(x+1)$, $f'(0) = \frac{d}{dx}f(x) \big|_{x=0}$. By \eqref{04}, we easily get
\begin{equation}\begin{split}\label{05}
\int_{\mathbb{Z}_p} e^{(x+y)t}   d\mu_{0} (y) = \frac{t}{e^t-1}e^{xt} = \sum_{n=0}^\infty B_n(x) \frac{t^n}{n!}.
\end{split}\end{equation}
Thus, from \eqref{05}, we have
\begin{equation}\begin{split}\label{06}
\int_{\mathbb{Z}_p} (x+y)^n   d\mu_{0} (y) = B_n(x),\quad (n \geq 0).
\end{split}\end{equation}
For $r \in \mathbb{N}$, the higher-order Bernoulli polynomials are defined by the multivariate $p$-adic invariant integrals on $\mathbb{Z}_p$ as follows:
\begin{equation}\begin{split}\label{07}
&\int_{\mathbb{Z}_p} \cdots \int_{\mathbb{Z}_p} e^{(x_1+x_2+\cdots+x_r+x)t}   d\mu_{0} (x_1) \cdots  d\mu_{0} (x_r) = \left( \frac{t}{e^t-1} \right)^r e^{xt} \\
&= \sum_{n=0}^\infty B_n^{(r)} (x) \frac{t^n}{n!}, \quad (\textnormal{see} \,\, [7,13]).
\end{split}\end{equation}
Comparing the coefficients on the both sides of \eqref{07}, we have
\begin{equation}\begin{split}\label{08}
\int_{\mathbb{Z}_p} \cdots \int_{\mathbb{Z}_p} (x_1+x_2+\cdots+x_r+x)^n   d\mu_{0} (x_1) \cdots  d\mu_{0} (x_r) = B_n^{(r)}(x),
\end{split}\end{equation}
where $B_n^{(r)}(x)$ are called higher-order Bernoulli polynomials.
In [2], L. Carlitz considered $q$-Bernoulli numbers which are given by
\begin{equation}\begin{split}\label{10}
\beta_{0,q}=1, \quad q(q\beta_q +1)^n - \beta_{n,q} = \begin{cases}1,&\text{if}\,\, n=1\\
0,&\text{if}\,\,n>1, \end{cases}
\end{split}\end{equation}
with the usual convention about replacing $\beta_q^n$ by $\beta_{n,q}$. Note that $\lim_{q \rightarrow 1} \beta_{n,q} = B_n, \,\,(n \geq 0).$ He also defined $q$-Bernoulli polynomials as follows:
\begin{equation}\begin{split}\label{11}
\beta_{n,q}(x) = \sum_{l=0}^n {n \choose l} [x]_q^{n-l} q^{lx} \beta_{l,q},\quad (\textnormal{see} \,\, [2,13]).
\end{split}\end{equation}
For $f \in UD(\mathbb{Z}_p)$, the $p$-adic $q$-integral on $\mathbb{Z}_p$ is defined by Kim to be
\begin{equation}\begin{split}\label{12}
\int_{\mathbb{Z}_p}  f(x)  d\mu_{q} (x)= \lim_{N \rightarrow \infty} \frac{1}{[p^N]_q} \sum_{x=0}^{p^N-1} f(x)q^x, \quad (\textnormal{see} \,\, [13]).
\end{split}\end{equation}
In [13], Kim proved that Carlitz's $q$-Bernoulli polynomials can be represented by the $p$-adic $q$-integral on $\mathbb{Z}_p$:
\begin{equation}\begin{split}\label{13}
\int_{\mathbb{Z}_p}  e^{[x+y]_q t}  d\mu_{q} (y) = \sum_{n=0}^\infty \beta_{n,q}(x) \frac{t^n}{n!}.
\end{split}\end{equation}
L. Carlitz considered the degenerate Bernoulli polynomials which are defined by the generating function to be
\begin{equation}\begin{split}\label{14}
\frac{t}{(1+\lambda t)^{\frac{1}{\lambda}}-1}(1+\lambda t)^{\frac{x}{\lambda}} = \sum_{n=0}^\infty B_{n,\lambda}^* (x)\frac{t^n}{n!},\quad (\textnormal{see} \,\, [3]).
\end{split}\end{equation}
When $x=0$, $B_{n,\lambda}^* = B_{n,\lambda}^* (0)$ are called degenerate Bernoulli numbers. Note that $\lim_{\lambda \rightarrow 0} B_{n,\lambda}^* (x) = B_n(x), \,\, (n \geq 0)$.

Recently, Kim defined (fully) degenerate Bernoulli polynomials which are different from Carlitz's degenerate Bernoulli polynomials:
\begin{equation}\begin{split}\label{15}
\int_{\mathbb{Z}_p} (1+\lambda t)^{\frac{x+y}{\lambda}}   d\mu_{0} (y) =& \frac{\log(1+\lambda t)^{\frac{1}{\lambda}}}{(1+\lambda t)^{\frac{1}{\lambda}} -1} (1+\lambda t)^{\frac{x}{\lambda}}\\
=& \sum_{n=0}^\infty B_{n,\lambda}(x) \frac{t^n}{n!}\quad (\textnormal{see} \,\, [17]).
\end{split}\end{equation}
Note that $\lim_{\lambda \rightarrow 0} B_{n,\lambda}(x) = B_n(x), \,\, (n \geq 0).$ Indeed,
\begin{equation}\begin{split}\label{16}
\sum_{n=0}^\infty B_{n,\lambda}^* (x) \frac{t^n}{n!} &= \frac{t}{\log(1+\lambda t)^{\frac{1}{\lambda}}} \int_{\mathbb{Z}_p} (1+\lambda t)^{\frac{x+y}{\lambda}}   d\mu_{0} (y) \\
&=\left( \sum_{l=0}^\infty b_l \frac{\lambda^l}{l!} t^l \right) \left( \sum_{m=0}^\infty B_{m,\lambda}(x) \frac{t^m}{m!} \right)\\
&= \sum_{n=0}^\infty \left( \sum_{l=0}^n {n \choose l} \lambda^l b_l B_{n-l,\lambda}(x) \right) \frac{t^n}{n!},
\end{split}\end{equation}
where $b_n$ is the $n$-th Bernoulli numbers of the second kind. Thus, by \eqref{16}, we get
\begin{equation}\begin{split}\label{17}
B_{n,\lambda}^* (x) = \sum_{l=0}^n {n \choose l} \lambda^l b_l B_{n-l,\lambda}(x), \,\, (n \geq 0).
\end{split}\end{equation}
The higher-order Carlitz's $q$-Bernoulli polynomials are given by the multivariate $p$-adic $q$-integral on $\mathbb{Z}_p$ as follows:
\begin{equation}\begin{split}\label{18}
&\int_{\mathbb{Z}_p} \cdots \int_{\mathbb{Z}_p} e^{[x_1+x_2+\cdots+x_r+x]_q t}   d\mu_{q} (x_1)\cdots  d\mu_{q} (x_r)\\
&= \sum_{n=0}^\infty \beta_{n,q}^{(r)} (x) \frac{t^n}{n!},\quad (\textnormal{see} \,\, [13]).
\end{split}\end{equation}
Note that $\lim_{q \rightarrow 1} \beta_{n,q}^{(r)} (x) = B_n^{(r)} (x), \,\, (n \geq 0)$.
Recently, Kim introduced higher-order degenerate $q$-Bernoulli polynomials which are derived from the $p$-adic $q$-integral on $\mathbb{Z}_p$ as follows:
\begin{equation}\begin{split}\label{19}
&\int_{\mathbb{Z}_p} \cdots \int_{\mathbb{Z}_p} (1+\lambda t)^{\frac{1}{\lambda}[x_1+x_2+\cdots+x_r+x]_q}  d\mu_{q} (x_1)\cdots  d\mu_{q} (x_r)\\
&= \sum_{n=0}^\infty \beta_{n,\lambda,q}^{(r)}(x) \frac{t^n}{n!},\quad (\textnormal{see} \,\, [17]).
\end{split}\end{equation}
Note that $\lim_{\lambda \rightarrow 0} \beta_{n,\lambda, q}^{(r)}(x) = \beta_{n,q}^{(r)}(x)$ and $\lim_{q \rightarrow 1} \beta_{n,\lambda,q}^{(r)}(x) = B_{n,\lambda}^{(r)}(x)$, where $B_{n,\lambda}^{(r)}(x)$ are the higher-order degenerate Bernoulli polynomials which are given by the generating function to be
\begin{equation}\begin{split}\label{20}
\left( \frac{\log(1+\lambda t)^{\frac{1}{\lambda}}}{(1+\lambda t)^{\frac{1}{\lambda}} -1} \right)^r (1+\lambda t)^{\frac{x}{\lambda}} = \sum_{n=0}^\infty B_{n,\lambda}^{(r)} (x) \frac{t^n}{n!},\quad (\textnormal{see} \,\, [3]).
\end{split}\end{equation}
The purpose of this paper is to give some identities of symmetry for the higher-order degenerate $q$-Bernoulli polynomials which are derived from the integral equations of the $p$-adic $q$-integral on $\mathbb{Z}_p$.

\section{Symmetric identities for the higher-order degenerate $q$-Bernoulli polynomials}

Let $\lambda, t \in \mathbb{C}_p$ with $|\lambda|_p \leq 1$ and $|t|_p < p^{-\frac{1}{p-1}}$ and let $w_1,w_2 \in \mathbb{N}$. From \eqref{04} and \eqref{20}, we note that
\begin{equation}\begin{split}\label{21}
\int_{\mathbb{Z}_p} \cdots \int_{\mathbb{Z}_p} (1+\lambda t)^{\frac{x_1+\cdots+x_r+x}{\lambda} }  d\mu_{0} (x_1)\cdots  d\mu_{0} (x_r)&= \left( \frac{\log(1+\lambda t)^{\frac{1}{\lambda}}}{(1+\lambda t)^{\frac{1}{\lambda}} -1} \right)^r (1+\lambda t)^{\frac{x}{\lambda}}
\\&= \sum_{n=0}^\infty B_{n,\lambda}^{(r)}(x) \frac{t^n}{n!}.
\end{split}\end{equation}
In view of \eqref{18}, we consider the higher-order degenerate $q$-Bernoulli polynomials which are derived from the multivariate $p$-adic $q$-integral on $\mathbb{Z}_p$ as follows:
\begin{equation}\begin{split}\label{22}
\int_{\mathbb{Z}_p} \cdots \int_{\mathbb{Z}_p} (1+\lambda t)^{\frac{1}{\lambda}[x_1+x_2+\cdots+x_r+x]_q}  d\mu_{q} (x_1)\cdots  d\mu_{q} (x_r)= \sum_{n=0}^\infty \beta_{n,\lambda,q}^{(r)}(x) \frac{t^n}{n!}.
\end{split}\end{equation}
Thus, we note that $\lim_{q \rightarrow 1} \beta_{n,\lambda,q}^{(r)}(x) = B_{n,\lambda}^{(r)}(x)$, $(n \geq 0)$. Now, we observe that
\begin{equation}\begin{split}\label{23}
&(1+\lambda t)^{\frac{1}{\lambda}[x_1+x_2+\cdots+x_r+x]_q} = \sum_{n=0}^\infty {\frac{1}{\lambda}[x_1+\cdots+x_r+x]_q \choose n} \lambda^n t^n\\
&=\sum_{n=0}^\infty \Big( \frac{1}{\lambda} [x_1+\cdots+x_r+x]_q \Big)_n \frac{\lambda^n t^n}{n!}\\
&= \sum_{n=0}^\infty [x_1+x_2+\cdots+x_r+x]_q \cdot \Big( [x_1+\cdots+x_r+x]_q - \lambda \Big) \cdots \Big( [x_1+\cdots+x_r+x]_q-(n-1)\lambda \Big) \frac{t^n}{n!}\\
&= \sum_{n=0}^\infty \Big( [x_1+\cdots+x_r+x]_q \Big)_{n,\lambda} \frac{t^n}{n!},
\end{split}\end{equation}
where $\Big([x]_q \Big)_{0,\lambda}=1$, $\Big([x]_q \Big)_{n,\lambda} = [x]_q \Big( [x]_q -\lambda \Big) \cdots \Big( [x]_q -(n-1)\lambda \Big)$, $(n \geq 1)$. It is not difficult to show that
\begin{equation}\begin{split}\label{24}
\Big( [x_1+\cdots+x_r+x]_q \Big)_{n,\lambda} = \sum_{l=0}^n S_1(n,l) \lambda^{n-l} [x_1+\cdots+x_r+x]_q^l,
\end{split}\end{equation}
where $S_1(n,l)$ is the stirling number of the first kind. By \eqref{18} and \eqref{22}, we get
\begin{equation}\begin{split}\label{25}
\beta_{n,\lambda,q}^{(r)} (x) &= \int_{\mathbb{Z}_p} \cdots \int_{\mathbb{Z}_p} \Big( [x_1+\cdots+x_r+x]_q \Big)_{n,\lambda} d\mu_{q} (x_1)\cdots  d\mu_{q} (x_r)\\
&= \sum_{l=0}^n S_1(n,l) \lambda^{n-l} \int_{\mathbb{Z}_p} \cdots \int_{\mathbb{Z}_p} [x_1+\cdots+x_r+x]_q^l  d\mu_{q} (x_1)\cdots  d\mu_{q} (x_r)\\
&=\sum_{l=0}^n S_1(n,l) \lambda^{n-l} \beta_{l,q}^{(r)} (x).
\end{split}\end{equation}
Let $w_1, w_2 \in \mathbb{N}$. Then, we have
\begin{equation}\begin{split}\label{26}
&\frac{1}{[w_1]_q^r} \int_{\mathbb{Z}_p} \cdots \int_{\mathbb{Z}_p} (1+\lambda t)^{
\frac{[w_1w_2x + w_2 \sum_{l=1}^r j_l + w_1 \sum_{l=1}^r y_l]_q}{\lambda}}
d\mu_{q^{w_1}} (y_1)\cdots  d\mu_{q^{w_1}} (y_r)\\
&=\frac{1}{[w_1]_q^r} \lim_{N \rightarrow \infty} \sum_{y_1,\cdots ,y_r =0}^{p^N-1} \frac{1}{[p^N]_{q^{w_1}}^r }(1+\lambda t)^{
\frac{[w_1w_2x + w_2 \sum_{l=1}^r j_l + w_1 \sum_{l=1}^r y_l]_q}{\lambda}} q^{w_1(y_1+\cdots y_r)}\\
&=\frac{1}{[w_1]_q^r} \lim_{N \rightarrow \infty} \frac{1}{[w_2 p^N]_{q^{w_1}}^r}\sum_{y_1,\cdots ,y_r =0}^{w_2p^N-1}(1+\lambda t)^{
\frac{[w_1w_2x + w_2 \sum_{l=1}^r j_l + w_1 \sum_{l=1}^r y_l]_q}{\lambda}} q^{w_1y_1 + w_1y_2 +\cdots +w_1y_r}\\
&=\lim_{N \rightarrow \infty} \frac{1}{[w_1w_2p^N]_q^r} \sum_{i_1,\cdots,i_r=0}^{w_2-1} \sum_{y_1,\cdots ,y_r =0}^{p^N-1} (1+\lambda t)^{
\frac{[w_1w_2x + w_2 \sum_{l=1}^r j_l + w_1 \sum_{l=1}^r (i_l+w_2 y_l)]_q}{\lambda}}\\
&\quad \times  q^{w_1(i_1+w_2y_1)+w_1(i_2+w_2y_2)+\cdots+w_1(i_r+w_2y_r)}\\
&=  \sum_{i_1,\cdots,i_r=0}^{w_2-1} q^{w_1 \sum_{l=1}^r i_l}\lim_{N \rightarrow \infty} \frac{1}{[w_1w_2p^N]_q^r}\sum_{y_1,\cdots ,y_r =0}^{p^N-1} (1+\lambda t)^{
\frac{[w_1w_2x + w_2 \sum_{l=1}^r j_l + w_1 \sum_{l=1}^r (i_l+w_2 y_l)]_q}{\lambda}}\\
&\quad \times q^{w_1w_2y_1 + w_1w_2y_2 +\cdots + w_1w_2y_r}.
\end{split}\end{equation}
By \eqref{26}, we get
\begin{equation}\begin{split}\label{27}
&I(w_1,w_2) = \frac{1}{[w_1]_q^r}  \sum_{i_1,\cdots,i_r=0}^{w_1-1} q^{w_2 \sum_{l=1}^r i_l}\int_{\mathbb{Z}_p} \cdots \int_{\mathbb{Z}_p} (1+\lambda t)^{
\frac{[w_1w_2x + w_2 \sum_{l=1}^r j_l + w_1 \sum_{l=1}^r y_l]_q}{\lambda}}
d\mu_{q^{w_1}} (y_1)\cdots  d\mu_{q^{w_r}} (y_r)\\
&=\lim_{N \rightarrow \infty} \frac{1}{[w_1w_2p^N]_q^r}\sum_{i_1,\cdots,i_r=0}^{w_2-1}
\sum_{j_1,\cdots,j_r=0}^{w_1-1}\sum_{y_1,\cdots ,y_r =0}^{p^N-1} q^{w_1 \sum_{l=1}^r i_l + w_2 \sum_{l=1}^r j_l + w_1w_2 \sum_{l=1}^r y_l}\\
&\quad \times (1+\lambda t)^{\frac{1}{\lambda}[w_1w_2x + w_2 \sum_{l=1}^r j_l + w_1 \sum_{l=1}^r i_l + w_1w_2 \sum_{l=1}^r y_l]_q  }.
\end{split}\end{equation}
On the other hand,
\begin{equation}\begin{split}\label{28}
&I(w_2, w_1) = \lim_{N \rightarrow \infty} \frac{1}{[w_1w_2p^N]_q^r}\sum_{i_1,\cdots,i_r=0}^{w_1-1} \sum_{j_1,\cdots,j_r=0}^{w_2-1}\sum_{y_1,\cdots ,y_r =0}^{p^N-1} q^{w_1 \sum_{l=1}^r j_l + w_2 \sum_{l=1}^r i_l + w_1w_2 \sum_{l=1}^r y_l}\\
&\quad \times (1+\lambda t)^{\frac{1}{\lambda}[w_1w_2x + w_1 \sum_{l=1}^r j_l + w_2 \sum_{l=1}^r i_l + w_1w_2 \sum_{l=1}^r y_l]_q  }.
\end{split}\end{equation}
From \eqref{27} and \eqref{28}, we note that $I(w_1,w_2) = I(w_2,w_1)$. Therefore, we obtain the following theorem.\\\\
\begin{theorem}
For $w_1,w_2 \in \mathbb{N}$, we have
\begin{equation*}\begin{split}
&\frac{1}{[w_2]_q^r} \sum_{j_1,\cdots,j_r=0}^{w_2-1} q^{w_1 \sum_{l=1}^r j_l}  \int_{\mathbb{Z}_p} \cdots \int_{\mathbb{Z}_p} (1+\lambda t)^{
\frac{1}{\lambda}[w_1w_2x + w_1 \sum_{l=1}^r j_l + w_2 \sum_{l=1}^r y_l]_q}
d\mu_{q^{w_2}} (y_1)\cdots  d\mu_{q^{w_2}} (y_r)
\\&=\frac{1}{[w_1]_q^r} \sum_{j_1,\cdots,j_r=0}^{w_1-1} q^{w_2 \sum_{l=1}^r j_l}  \int_{\mathbb{Z}_p} \cdots \int_{\mathbb{Z}_p} (1+\lambda t)^{
\frac{1}{\lambda}[w_1w_2x + w_2 \sum_{l=1}^r j_l + w_1 \sum_{l=1}^r y_l]_q}
d\mu_{q^{w_1}} (y_1)\cdots  d\mu_{q^{w_1}} (y_r).
\end{split}\end{equation*}
\end{theorem}
It is not difficult to show that
\begin{equation}\begin{split}\label{29}
\bigg[w_1w_2x + \sum_{l=1}^r j_l w_2 + \sum_{l=1}^r y_l w_1\bigg]_q  = [w_1]_q \bigg[w_2x + \frac{w_2}{w_1} \sum_{l=1}^r j_l + \sum_{l=1}^r y_l\bigg]_{q^{w_1}}.
\end{split}\end{equation}
Thus, by \eqref{22} and \eqref{29}, we get
\begin{equation}\begin{split}\label{30}
&\int_{\mathbb{Z}_p} \cdots \int_{\mathbb{Z}_p} (1+\lambda t)^{\frac{1}{\lambda}[w_1w_2x + \sum_{l=1}^r j_l w_2 + \sum_{l=1}^r y_l w_1]_q}  d\mu_{q^{w_1}} (y_1)\cdots  d\mu_{q^{w_1}} (y_r)\\
&=\int_{\mathbb{Z}_p} \cdots \int_{\mathbb{Z}_p} (1+\lambda t)^{\frac{[w_1]_q}{\lambda}
[w_2x + \frac{w_2}{w_1} \sum_{l=1}^r j_l + \sum_{l=1}^r y_l]_{q^{w_1}} }  d\mu_{q^{w_1}} (y_1)\cdots  d\mu_{q^{w_1}} (y_r)\\
&=\int_{\mathbb{Z}_p} \cdots \int_{\mathbb{Z}_p} \left( 1+ \frac{\lambda}{[w_1]_q}[w_1]_q t \right)^{\frac{[w_1]_q}{\lambda} [w_2x + \frac{w_2}{w_1} \sum_{l=1}^r j_l + \sum_{l=1}^r y_l]_{q^{w_1}} }  d\mu_{q^{w_1}} (y_1)\cdots  d\mu_{q^{w_1}} (y_r)\\
&= \sum_{n=0}^\infty \beta_{n,\frac{\lambda}{[w_1]_q},q^{w_1}}^{(r)} \Big(w_2x + \frac{w_2}{w_1} \sum_{l=1}^r j_l\Big) [w_1]_q^n \frac{t^n}{n!}.
\end{split}\end{equation}
Therefore, by Theorem 1 and \eqref{30}, we obtain the following theorem.\\\\
\begin{theorem}
Let $n \geq 0$ and $w_1,w_2 \in \mathbb{N}$. Then we have
\begin{equation*}\begin{split}
&[w_1]_q^{n-r} \sum_{j_1,\cdots,j_r=0}^{w_1-1} q^{w_2 \sum_{l=1}^r j_l} \beta_{n,\frac{\lambda}{[w_1]_q},q^{w_1}}^{(r)}\Big( w_2x + \frac{w_2}{w_1} \sum_{l=1}^r j_l \Big)
\\&=[w_2]_q^{n-r} \sum_{j_1,\cdots,j_r=0}^{w_2-1} q^{w_1 \sum_{l=1}^r j_l} \beta_{n,\frac{\lambda}{[w_2]_q},q^{w_2}}^{(r)}\Big( w_1x + \frac{w_1}{w_2} \sum_{l=1}^r j_l \Big).
\end{split}\end{equation*}
\end{theorem}
In particular, if we take $w_2=1$, then we have
\begin{equation}\begin{split}\label{31}
[w_1]_q^{n-r} \sum_{j_1,\cdots,j_r=0}^{w_1-1} q^{\sum_{l=1}^r j_l} \beta_{n,\frac{\lambda}{[w_1]_q},q^{w_1}}^{(r)}\Big(x + \frac{1}{w_1} \sum_{l=1}^r j_l \Big) = \beta_{n,\lambda,q}^{(r)} (w_1x).
\end{split}\end{equation}
From \eqref{31}, we note that
\begin{equation*}\begin{split}
 \beta_{n,q}^{(r)} (w_1x) &= \lim_{\lambda \rightarrow 0} \beta_{n,\lambda,q}^{(r)} (w_1x) = \lim_{\lambda \rightarrow 0} [w_1]_q^{n-r} \sum_{j_1,\cdots,j_r=0}^{w_1-1} q^{\sum_{l=1}^r j_l} \beta_{n,\frac{\lambda}{[w_1]_q},q^{w_1}}^{(r)}\Big(x + \frac{1}{w_1} \sum_{l=1}^r j_l \Big)\\
 &=[w_1]_q^{n-r} \sum_{j_1,\cdots,j_r=0}^{w_1-1} q^{\sum_{l=1}^r j_l} \beta_{n,q^{w_1}}^{(r)}\Big(x + \frac{1}{w_1} \sum_{l=1}^r j_l \Big).
\end{split}\end{equation*}
From \eqref{25}, we note that
\begin{equation}\begin{split}\label{32}
&\int_{\mathbb{Z}_p} \cdots \int_{\mathbb{Z}_p} \left([w_1]_q \bigg[w_2x + \frac{w_2}{w_1} \sum_{l=1}^r j_l + \sum_{l=1}^r y_l\bigg]_{q^{w_1}}\right)_{n,\lambda} d\mu_{q^{w_1}} (y_1)\cdots  d\mu_{q^{w_1}} (y_r)\\
&=\sum_{l=0}^n S_1(n,l) [w_1]_q^l \lambda^{n-l}\int_{\mathbb{Z}_p} \cdots \int_{\mathbb{Z}_p}  \bigg[w_2x + \frac{w_2}{w_1} \sum_{l=1}^r j_l + \sum_{l=1}^r y_l\bigg]_{q^{w_1}}^l d\mu_{q^{w_1}} (y_1)\cdots  d\mu_{q^{w_1}} (y_r)\\
&= \sum_{l=0}^n S_1(n,l) [w_1]_q^l \lambda^{n-l} \sum_{i=0}^l {l \choose i} \left( \frac{[w_2]_q}{[w_1]_q} \right)^i [j_1+\cdots+j_r]_{q^{w_2}}^i q^{w_2(l-i) \sum_{k=1}^r j_k}\\
&\quad \times\int_{\mathbb{Z}_p} \cdots \int_{\mathbb{Z}_p} [w_2x + y_1+ \cdots y_r]_{q^{w_1}}^{l-i}  d\mu_{q^{w_1}} (y_1)\cdots  d\mu_{q^{w_1}} (y_r)\\
&=\sum_{l=0}^n \sum_{i=0}^l S_1(n,l)\lambda^{n-l} [w_1]_q^{l-i} [w_2]_q^i  {l \choose i} [j_1+\cdots+j_r]_{q^{w_2}}^i q^{w_2(l-i) \sum_{k=1}^r j_k} \beta_{l-i,q^{w_1}}^{(r)} (w_2x).
\end{split}\end{equation}
By \eqref{32}, we get
\begin{equation}\begin{split}\label{33}
&[w_1]_q^{n-r} \sum_{j_1,\cdots,j_r=0}^{w_1-1} q^{w_2 \sum_{l=1}^r j_l} \int_{\mathbb{Z}_p} \cdots \int_{\mathbb{Z}_p} \left([w_1]_q \bigg[w_2x + \frac{w_2}{w_1} \sum_{l=1}^r j_l + \sum_{l=1}^r y_l\bigg]_{q^{w_1}}\right)_{n,\lambda} d\mu_{q^{w_1}} (y_1)\cdots  d\mu_{q^{w_1}} (y_r)\\
&=\sum_{l=0}^n \sum_{i=0}^l {l \choose i} S_1(n,l)\lambda^{n-l} [w_1]_q^{n+l-r-i} [w_2]_q^i \beta_{l-i,q^{w_1}}^{(r)}(w_2x) \sum_{j_1,\cdots,j_r=0}^{w_1-1} [j_1+\cdots+j_r]_{q^{w_2}}^i q^{w_2(l+1-i) \sum_{k=1}^r j_k}  \\
&=\sum_{l=0}^n \sum_{i=0}^l {l \choose i} S_1(n,l)\lambda^{n-l} [w_1]_q^{n+l-r-i} [w_2]_q^i \beta_{l-i,q^{w_1}}^{(r)}(w_2x) T_{l+1,i}^{(r)} (w_1| q^{w_2})
\end{split}\end{equation}
where
\begin{equation}\begin{split}\label{34}
T_{n,i}^{(r)} (w|q) = \sum_{j_1, \cdots, j_r=0}^{w-1} q^{(n-i) \sum_{l=1}^r j_l} [j_1+\cdots+j_r]_q^i .
\end{split}\end{equation}
Therefore, by Theorem 2 and \eqref{34}, we obtain the following theorem.\\\\
\begin{theorem}
For $n \geq 0$, $w_1, w_2 \in \mathbb{N}$, we have
\begin{equation*}\begin{split}
&\sum_{l=0}^n \sum_{i=0}^l {l \choose i} S_1(n,l)\lambda^{n-l} [w_1]_q^{n+l-r-i} [w_2]_q^i \beta_{l-i,q^{w_1}}^{(r)}(w_2x) T_{l+1,i}^{(r)} (w_1| q^{w_2})\\
&=\sum_{l=0}^n \sum_{i=0}^l {l \choose i} S_1(n,l)\lambda^{n-l} [w_2]_q^{n+l-r-i} [w_1]_q^i \beta_{l-i,q^{w_2}}^{(r)}(w_1x) T_{l+1,i}^{(r)} (w_2| q^{w_1}).
\end{split}\end{equation*}
\end{theorem}
Let us take $\lambda \rightarrow 0$. Then, we have
\begin{equation*}\begin{split}
&\sum_{i=0}^n {n \choose i} [w_1]_q^{2n-r-i} [w_2]_q^i \beta_{n-i,q^{w_1}}^{(r)}(w_2x) T_{n+1,i}^{(r)} (w_1| q^{w_2})\\
&=\sum_{i=0}^n {n \choose i} [w_2]_q^{2n-r-i} [w_1]_q^i \beta_{n-i,q^{w_2}}^{(r)}(w_1x) T_{n+1,i}^{(r)} (w_2| q^{w_1}).
\end{split}\end{equation*}

\noindent
{\bf{Remark.}} Recently, several authors have studied the symmetry identities of Bernoulli and Euler polynomials and Carlitz's $q$-Bernoulli and $q$-Euler polynomials(see [1-24]).

\bibliographystyle{amsplain}

\begin{thebibliography}{10}

\bibitem{key-01} A. Adelberg, \textit{A finite difference approach to degenerate Bernoulli and Stirling polynomials}, Discrete Math., {\bf{140}} (1995), no. 1-3, 1-21.

\bibitem{key-02} L. Carlitz, \textit{$q$-Bernoulli and Eulerian numbers
}, Trans. Amer. Math. Soc., {\bf{76}}  (1954), 332--350.

\bibitem{key-03} L. Carlitz, \textit{Degenerate Stirling, Bernoulli and Eulerian numbers}, Utilitas Math., {\bf{15}} (1979), 51--88.

\bibitem{key-04} D. V. Dolgy, T. Kim, H.-I. Kwon, J. J. Seo,
\textit{Symmetric identities of degenerate $q$-Bernoulli polynomials under symmetry group $S_3$}, Proc. Jangjeon Math. Soc., {\bf{19}} (2016), no. 1, 1--9.

\bibitem{key-05} Y. He,  \textit{Symmetric identities for Carlitz's $q$-Bernoulli numbers and polynomials}, Adv. Difference Equ., {\bf{2013}} 2013:246, 10 pages.

\bibitem{key-06} F. T. Howard, \textit{Explicit formulas for degenerate Bernoulli numbers}, Discrete Math., {\bf{162}} (1996), no. 1-3, 175--185.

\bibitem{key-07} D. S. Kim, N. Lee, J. Na, K. H. Park,
\textit{Abundant symmetry for higher-order Bernoulli polynomials (I)
}, Adv. Stud. Contemp. Math., {\bf{23}} (2013), no. 3, 461--482.

\bibitem{key-08} D. S. Kim, N. Lee, J. Na, K. H. Park,
\textit{Identities of symmetry for higher-order Euler polynomials in three variables (I)
}, Adv. Stud Contemp. Math., {\bf{22}} (2012), no. 1, 51--74.

\bibitem{key-09} D. S. Kim,
\textit{Symmetry identities for generalized twisted Euler polynomials twisted by unramified roots of unity}, Proc. Jangjeon Math. Soc., {\bf{15}} (2012), no. 3, 303--316.

\bibitem{key-10} D. S. Kim, T. Kim,
\textit{Some identities of symmetry for $q$-Bernoulli polynomials under symmetric group of degree $n$}, Ars Comb., {\bf{126}} (2016), 435--441.

\bibitem{key-11} D. S. Kim, T. Kim,
\textit{Barnes-type Boole polynomials}, Contrib. Discrete Math., {\bf{11}} (2016), no. 1, 7--15.

\bibitem{key-12} D. S. Kim, N. Lee, J. Na, K. H. Park,
\textit{Abundant symmetry for higher-order Bernoulli polynomials (II)
}, Proc. Jangjeon Math. Soc., {\bf{16}} (2013), no. 3, 359--378.

\bibitem{key-13} T. Kim, \textit{$q$-Volkenborn integration}, Russ. J. Math. Phys., {\bf{9}} (2002), no. 3, 288-299.
\
\bibitem{key-14} T. Kim, D. S. Kim,
\textit{Identities involving degenerate Euler numbers and polynomials arising from non-linear differential equations}, J. Nonlinear Sci. Appl., {\bf{9}} (2016), no. 5, 2086--2098.

\bibitem{key-15} T. Kim, D. S. Kim,
\textit{Identities of symmetry for $(h,q)$- extensions of generalized higher-order Euler polynomials}, Adv. Stud. Contemp. Math., {\bf{26}} (2016), no. 3, 579--585.

\bibitem{key-16} T. Kim, J. J. Seo,
\textit{Revisited nonlinear differential equations arising from the generating funcions of degenerate Bernoulli numbers}, Adv. Stud. Contemp. Math., {\bf{26}} (2016), no. 3, 401--406.

\bibitem{key-17} T. Kim, \textit{On degenerate $q$-Bernoulli polynomials}, Bull. Korean Math. Soc. {\bf{53}} (2016), no. 4, 1149-1156.

\bibitem{key-18} T. Kim,
\textit{Symmetry of power sum polynomials and multivariate fermionic $p$-adic invariant integral on $\mathbb{Z}_p$}, Russ. J. Math. Phys., {\bf{16}} (2009), no. 1, 93--96.

\bibitem{key-19} Y.-H. Kim, K.-W. Hwang,
\textit{Symmetry of power sum and twisted Bernoulli polynomials},
Adv. Stud. Contemp. Math., {\bf{18}} (2009), no. 2, 127--133.

\bibitem{key-20} G. D. Liu,
\textit{Degenerate Bernoulli numbers and polynomials of higher order
}, (Chinese), J. Math. (Wuhan), {\bf{25}} (2005), no. 3, 283--288.

\bibitem{key-21} E.-J. Moon, S.-H. Rim, J.-H. Jin, S.-J. Lee,
\textit{On the symmetric properties of higher-order twisted $q$-Euler numbers polynomials},Adv. Difference Equ., {\bf{2010}} Art. ID 765259, 8 pages.


\bibitem{key-22} C. S. Ryoo,
\textit{A note on the weighted $q$-Euler numbers and polynomials
},Adv. Stud. Contemp. Math., {\bf{21}} (2011), no. 1, 47--54.

\bibitem{key-23} H. M. Srivastava, T. Kim, Y. Simsek,
\textit{$q$-Bernoulli numbers and polynomials associated with multiple $q$-zeta functions and bases $L$-series
}, Russ. J. Math. Phys., {\bf{12}} (2005), no. 2, 241--268.

\bibitem{key-24} Z. Zhang, J. Yang,
\textit{On sums of products of the degenerate Bernoulli numbers
}, Integral Transforms Spec. Funct., {\bf{20}} (2009), no. 9-10, 751--755.


\end{thebibliography}
\providecommand{\bysame}{\leavevmode\hbox to3em{\hrulefill}\thinspace}
\providecommand{\MR}{\relax\ifhmode\unskip\space\fi MR }
\providecommand{\MRhref}[2]{%
  \href{http://www.ams.org/mathscinet-getitem?mr=#1}{#2}
}
\providecommand{\href}[2]{#2}

\end{document}